\documentclass[12pt]{article}

\newtheorem{lem}{Lemma}[section]
\newtheorem{theorem}[lem]{Theorem }
\newtheorem{remark}[lem]{Remark }
\newtheorem{example}[lem]{Example }

\newtheorem{fact}[lem]{Fact }
\newtheorem{notation}[lem]{Notation }
\newtheorem{questions}[lem]{Questions }
\newtheorem{corollary}[lem]{Corollary}

\usepackage{amsfonts} \usepackage{amssymb}
\newcommand{\UUU}{\mathcal{U}} \newcommand{\NNN}{\mathcal{N}}
\newcommand{\AAA}{\mathcal{A}}\newcommand{\BBB}{\mathcal{B}}\newcommand{\CCC}{\mathcal{C}}
\newcommand{\ifff}{\Longleftrightarrow}
\newcommand{\mod}{\mbox{mod~}}  \newcommand{\card}{\mbox{card~}}

\newcommand{\ZZ}{\mathbb{Z}} 
 \newcommand{\BB}{\mathbb{B}}

\begin{document}

\baselineskip 10 pt

\begin{center}
{\bf \large Conjugacy classes and invariant subrings \\ \vskip 8 pt
of $R$-automorphisms of $R[x]$} \footnote{This work is supported by
a research grant from Yarmouk University} \footnote{2000 AMS
Classification Number: 13A50}
\end{center}

\bigskip
\begin{center}
{\bf Jebrel M. Habeb and Mowaffaq Hajja}\\
{Yarmouk University}\\
{Irbid -- Jordan}\\
{jhabeb@yu.edu.jo~,~~mhajja@yu.edu.jo}
\end{center}

\begin{center} {and}
\end{center}

\begin{center}
{\bf William J. Heinzer}\\
{Purdue University}\\
{West Lafayette, IN 47907 -- USA}\\
{heinzer@math.purdue.edu}
\end{center}

\medskip
\baselineskip 16 pt

\section{Introduction and terminology}

All rings are assumed to be commutative with an identity element.
The group of units of a ring $R$ is denoted by $\UUU (R)$, and the
set of nilpotent elements of $R$ by $\NNN (R)$. It is well known
that $\NNN (R)$ is an ideal of $R$ called  the {\it nilradical} of
$R$. If $\NNN (R) = \{ 0 \}$, then  $R$ is said to be {\it reduced}.

Let  $R[x]$ be the polynomial ring in one indeterminate $x$ over a
ring $R$. An endomorphism $\sigma$ of $R[x]$ is called an
$R$-endomorphism if $\sigma(r)= r$ for all $r$ in $R$. Clearly, an
$R$-endomorphism $\sigma$ of $R[x]$ is completely determined by
$\sigma (x)$. A theorem of Gilmer \cite[Theorem 3]{Gilmer} asserts
the following:

\begin{fact} \label{1.0}
 An $R$-endomorphism $\sigma$ is an $R$-automorphism if
and only if
\begin{eqnarray} \label{Gilmer}
    \sigma (x) = a + u x + x^2 f(x),
\end{eqnarray}
where  $a \in R$, $u \in \UUU (R)$, and $f(x) \in \NNN(R[x])$. In
other words, an element  $y \in R[x]$ is such that $R[y] = R[x]$
if and only if $y$ is of the form given by the right hand side of
(\ref{Gilmer}).
\end{fact}

It is well known that a polynomial $f(x) \in R[x]$ is nilpotent if
and only if all the coefficients of $f$ are nilpotent elements in
$R$ (see \cite[Exercise 2(ii), page 11]{AM}). Thus  $\NNN(R[x]) =
\NNN(R)R[x]$.  Since the sum of a unit and a nilpotent element is a
unit,  an equivalent formulation of Fact \ref{1.0} is:

\begin{fact} \label{1.00}
An $R$-endomorphism $\sigma$ of $R[x]$ is an $R$-automorphism if and
only if
\begin{eqnarray} \label{Gilmer2}
    \sigma (x) = b + v x + g(x),
\end{eqnarray}
where $b \in R$, $v \in \UUU(R)$ and $g(x) \in \NNN(R[x])$.
\end{fact}

\begin{remark} \label{1.1}{\em
Let $R$ be an integral domain with field of fractions $K$  and let
$H$ be a group of $R$-automorphism of $R[x]$. Each $h \in H$ extends
in a canonical way to an automorphism of the field $K(x)$. Thus we
may regard $H$ as a group of automorphism of the field $K(x)$. The
fixed field $K(x)^H$ of $H$ acting on $K(x)$ contains the fixed ring
$R[x]^H$ of $H$ acting on $R[x]$. If $H$ is infinite, then the fixed
field $K(x)^H$ is $K$. Therefore if $H$ is an infinite group of
$R$-automorphisms of $R[x]$, where $R$ is an integral domain, then
$R$ is the ring of invariants of $H$ acting on $R[x]$, i.e., $R[x]^H
= R$. Assume the group $H$ is finite, say $|H| = n$, and let $L =
K(x)^H$. Then $K(x)/L$ is a Galois field extension with $[K(x):L] =
n$ and $\{1, x, \ldots,x^{n-1}\}$ is a vector space basis for $K(x)$
over $L$. Moreover, $L$ is the field of fractions of $R[x]^H$
\cite[Corollary, page 324]{Bourbaki}. For each $h \in H$, we have
$h(x) = u_hx + a_h$, where $u_h \in \UUU(R)$ and $a_h \in R$. Let $f
= \prod_{h \in H}(u_hx + a_h)$ denote the norm of $x$ with respect
to $H$. Notice that $$ f = ux^n + b_{n-1}x^{n-1} + \cdots + b_1x +
b_0,
$$
where $u \in \UUU(R)$ and each $b_i \in R$. It follows that $x$
satisfies a monic polynomial of degree $n$ with coefficients in
$R[f]$. Therefore $K(f) = L$ and $\{1,x, \ldots, x^{n-1}\}$ is a
free module basis for $R[x]$ as an $R[f]$-module, and as Samuel
observes in \cite{Samuel}, we must have $R[x]^H = R[f]$.

}
\end{remark}

For an  integer $n \ge 2$, let $\ZZ_n$ denote the residue class
ring $\ZZ/n\ZZ$. Notice that any automorphism of the polynomial
ring $\ZZ_n[x]$ maps $1$ to $1$ and thus maps every element
of $\ZZ_n$ to itself and is therefore a $\ZZ_n$-automorphism
of $\ZZ_n[x]$. . If $p$ is a prime integer, and $G$ is the
group of automorphism of $\ZZ_p[x]$, then the result of Samuel in
Remark \ref{1.1} implies that the ring of invariants $\ZZ_p[x]^G =
\ZZ_p[(x^p-x)^{p-1}]$

Let $H$ be a finite group of $R$-automorphisms of the polynomial
ring $R[x]$.  Example \ref{simple} illustrates that computing
generators over $R$ for the ring of invariants $R[x]^H$ is more
subtle in the case where $R$ has nonzero nilpotent elements. With $R
= \ZZ_4$, there exists a group $H = \langle \alpha \rangle$ of order 2 of
$R$-automorphisms of $R[x]$  such that $R[x]^H$ properly contains
the subring $R[x\alpha(x), x + \alpha(x)]$ generated over $R$ by the norm and
trace of $x$ with respect to $H$.

\begin{example} \label{simple} {\em
Let $\alpha$  denote the  automorphism of $\ZZ_4[x]$ defined by
$\alpha(x) = -x$ and let $H = \langle \alpha \rangle$ denote the cyclic
group generated by $\alpha$. Then $H$ has order 2 and $\ZZ_4[x + \alpha(x),
x\alpha(x)] = \ZZ_4[x^2]$ is properly contained in $\ZZ_4[x]^H $; for
it is clear that $\ZZ_4[2x, x^2] \subseteq \ZZ_4[x]^H$. Indeed,
$\ZZ_4[x]^H = \ZZ_4[2x, x^2]$, as we discuss in more detail below
in Section 3.

}
\end{example}

\begin{notation} \label{1.4} {\em
In \cite{Dowlen1}, the group of $R$-automorphisms of $R[x]$ is
denoted by $G(R)$, and the subgroup of $G(R)$ consisting of
$R$-automorphisms in which $f(x)$, as given in (\ref{Gilmer}), is
zero by $B(R)$. We retain the notation $G(R)$; however, to highlight
the dependence of $B(R)$ on $x$, we denote it by $B_x(R)$, or rather
by $\BB_x(R)$. This  dependence is illustrated in Example \ref{2.1},
where it is shown that if $R[x] = R[y]$, then $\BB_x(R)$ need not
coincide with $\BB_y(R)$. Of course it is true  that $\BB_x(R)$ and
$\BB_y(R)$ are isomorphic via the obvious  map that sends the
$R$-automorphism $s$ defined by $s(x) = a + ux$ to the
$R$-automorphism $s'$ defined by $s'(y) = a + uy$. Thus, up to
isomorphism, one may denote $\BB_x(R)$ by $\BB (R)$. Identifying the
element of $\BB (R)$ defined by $x \mapsto ux + a$ with the element
$(u,a) \in \UUU (R) \times R$,  we see that $\BB (R)$ is the
semidirect product of the multiplicative group $\UUU (R)$ by the
additive group $R$ defined by the multiplication
\begin{eqnarray*}
(u,a) \cdot (v,b) = (uv,va+b).
\end{eqnarray*}

}
\end{notation}

At this point, it is convenient to prove the following  simple
theorem that will be used later in the proof of Theorem \ref{T30:6}.

\begin{theorem} \label{products}
{\it Let  $R \times S$  denote the direct product of the rings $R$ and $S$.
Then the group $\BB(R \times S)$ is isomorphic to the direct product
$\BB(R)  \times \BB(S)$ of the groups $\BB(R)$ and $\BB(S)$.}
\end{theorem}

\noindent {\it Proof.} It is easy to check that the mapping from
$R[x] \times S[x]$ to $(R \times S) [x]$ defined by
$$\left( \sum r_i x^i , \sum s_i x^i \right) \mapsto \sum (r_i,s_i) x^i$$
is a ring isomorphism and that the mapping from
$\BB_x (R) \times \BB_x (S)$ to $\BB_x (R \times S)$ defined by
$$\left( x \mapsto ux + a , x \mapsto vx+b \right) \mapsto
\left(x \mapsto (u,v)x + (a,b) \right)$$
is a group isomorphism.
\hfill $\Box$

\bigskip

\begin{remark} \label{1.6} {\em
If the ring $R$ has nonzero nilpotent elements, then the
group $G = G(R)$ of automorphisms of $R[x]$ is infinite, and
properly contains the subgroup $\BB_x(R)$ defined in (\ref{1.4});
however, Dowlen proves that the  ring of invariants $R[x]^{G}$ is
equal to the ring of invariants $R[x]^{\BB_x}$ \cite[Theorem
1.2]{Dowlen1}. In the case where $R$ is a finite ring, for example
$R = \ZZ_n$, the group $\BB_x(R)$ is finite and $R[x]$ is integral
over its invariant subring $R[x]^G = R[x]^{\BB_x}$. Dowlen in
\cite{Dowlen2} determines generators for the  ring of invariants
of $\ZZ_n[x]$ with respect to the group $G(\ZZ_n)$ of
automorphisms of $\ZZ_n[x]$. One of our goals is to determine
generators for the ring of  invariants  of $\ZZ_n[x]$ with respect
to various subgroups $H$ of $G(\ZZ_n)$. }
\end{remark}

In the hope that it may be useful in the task of describing rings of
invariants of $R[x]$, we prove  in Section 4 that every element of
$\BB_x(\ZZ_n)$ is equivalent to an element having a certain simple
representation. We determine conditions for two elements of
$\BB_x(\ZZ_n)$ to be conjugate and give a formula for the number of
conjugacy classes of this group. In Section 2, we state and prove
results that hold in general, and raise several open problems
concerning $\BB_x(R)$ and $G(R)$. We examine in detail in Section 3
the structure of the automorphism group $G(\ZZ_4)$ of $\ZZ_4[x]$ and
enumerate the invariant subrings of $\ZZ_4[x]$ with respect to
subgroups of $G(\ZZ_4)$. In particular,  for $R = \ZZ_4$, we
establish the existence of subrings of $R[x]$ that are rings of
invariants of  subgroups of $G(R)$, but are not rings of invariants
of subgroups of $\BB_x(R)$. We prove, however, that each of these
invariant subrings of $R[x]$ is the ring of invariants of a subgroup
of $\BB_z(R)$ for some $z \in R[x]$ such that $R[z] = R[x]$.

\section{General results and open problems}

Let $A$ be the ring of polynomials in one indeterminate over the
ring $R$, and let $G = G(R)$ be the group of $R$-automorphisms of
$A$. For $x \in A$ such that $A=R[x]$ and for $\sigma \in G$, we say
that $\sigma$ is $x$-{\it basic} if $\sigma (x)$ is of the form
$\sigma (x) = a + ux$ for some $a \in R$ and $u \in \UUU (R)$. The
group of $x$-basic elements of $G$ is denoted by $\BB_x=\BB_x(R)$.
As mentioned in (\ref{1.4}), $\BB_x$ and $\BB_y$ are isomorphic if
$y$ is such that $R[x] = R[y]$. However, they need not be  equal as
is illustrated in Example \ref{2.1}.

\begin{example} \label{2.1} {\em
Let $n = p^k$, where $p \ge 3$ is an odd prime and where $k \ge
2$. Let $R = \ZZ_n$ and let  $\sigma \in G(R)$ be defined by
$\sigma (x) = 2x.$ Then $\sigma$ is $x$-basic. However, if we let
$y = x + p^{k-1}x^2$, then it follows from Fact \ref{1.0} that
$R[y] = R[x]$, and it is easy to see that $\sigma$ is not
$y$-basic. Indeed:
 \begin{eqnarray*}
    \sigma (y) &=& 2x + p^{k-1}(2x)^2 ~=~ 2y + 2 p^{k-1} x^2 ~=~
    2y + 2 p^{k-1}  (y-p^{k-1} x^2)^2\\&=&
        2y + 2p^{k-1} y^2.
\end{eqnarray*}

}
\end{example}

Two elements $\alpha$ and $\beta$ in $G(R)= G$ are $G$-conjugate
if $\alpha = g^{-1} \beta g$ for some $g \in G(R)$. They are
$\BB_x$-conjugate if $g$ can be chosen to belong to $\BB_x(R)$.
They are $\BB$-conjugate if $g$ can be chosen to belong to
$\BB_y(R)$ for some $y$ such that $R[x] = R[y]$. The corresponding
conjugacy classes are  denoted by $[\alpha]_G$,
$[\alpha]_{\BB_x}$, and $[\alpha]_{\BB}$, respectively.

If $\alpha \in \BB_x(R)$, then we clearly have $ [\alpha]_{\BB_x}
\subseteq [\alpha]_{\BB} \subseteq [\alpha]_G.$  Theorem \ref{2.22}
demonstrates  that if $R$ is not reduced and $\alpha \in \BB_x$ is
defined by $\alpha(x) = x+1$, then $[\alpha]_{\BB_x} \subsetneq
[\alpha]_G$.

\begin{theorem} \label{2.22}  {\it Let
$R$ be a ring that is not reduced and let $\alpha \in \BB_x$ be
defined by $\alpha(x) = x+1$. Then $[\alpha]_G$ is not contained
in $\BB_x$. Therefore $\BB_x(R)$ is not normal in $G(R)$.}
\end{theorem}

\noindent {\it Proof.} Since $R$ is not reduced,   there exists a
nonzero $r \in R$ such that $r^2 = 0$. If $[\frac{k(k-1)}{2}]r=0$
for 3 consecutive natural numbers $k$, say
$$\frac{n(n-1)}{2}~r = \frac{(n+1)n}{2}~r = \frac{(n+2)(n+1)}{2}~r = 0,$$
then  by subtracting we have $n r = (n+1)r =0$. But this implies
that $r = 0$. Therefore there exists an integer $n \ge 3$ such that
$[\frac{n(n-1)}{2}]r \ne 0$. For such  $r$ and $n$, define $\sigma$
in $G(R)$  by
\begin{eqnarray*}
    \sigma (x) &=& x + rx^2+r x^n.
\end{eqnarray*}
We  prove that $\sigma^{-1} \alpha \sigma \notin \BB_x(R)$. For
suppose that
\begin{eqnarray}
\sigma^{-1} \alpha  \sigma (x) = ux + c, \label{EQ1}
\end{eqnarray}
where $u \in \UUU(R)$ and $c \in R$.  From (\ref{EQ1}), it follows
that
\begin{eqnarray*}
    \alpha (x + rx^2+r x^n) &=& u(x+rx^2+rx^n) + c\\
    x+1 + r(x+1)^2+ r(x+1)^n &=& u(x+rx^2+rx^n)+c.
    \end{eqnarray*}
Equating the coefficients of $x$ and of $x^2$, we see that
$$1+2r+nr=u~~~~~\mbox{ and }~~~~~ r+r\frac{n(n-1)}{2}=ur.$$
Hence $$r+2r^2+nr^2 = r+r\frac{n(n-1)}{2},$$ and  therefore
$[\frac{n(n-1)}{2}] r = 0$, contradicting the choice of $n$. \hfill
$\Box$

\begin{remark} \label{2.3} {\em
Let $ \alpha, \beta \in \BB_x(R)$ be defined by
\begin{eqnarray}  \label{2.31}
    \alpha(x) = ux + a, \mbox{\hskip 15 pt} \beta(x) = vx + b,
\end{eqnarray}
where $u,v \in \UUU(R)$ and $a,b \in R$. It is easy to see that if
$\alpha$ and $\beta$ are  $\BB_x$-conjugate, then $u = v$. For if
$\sigma \in \BB_x(R)$ is defined by $\sigma(x) = wx + c$, then
$\sigma^{-1}(x) = w^{-1}x - w^{-1}c$,  and
\begin{eqnarray*}
  \sigma^{-1}\alpha\sigma(x)
    &=&  \sigma^{-1}\alpha(wx+c)   \\
&=& \sigma^{-1}(ws(x)+c)    \\
&=&  \sigma^{-1}(wux + wa + c)\\
&=& wu\sigma^{-1}(x) + wa + c\\
&=&  wuw^{-1}x - wuw^{-1}c +wc + c\\
&=& ux - uc + wc + c.
                      \end{eqnarray*}
Therefore if $\alpha$ and $\beta$ in $\BB_x(R)$  as in (\ref{2.31})
are $\BB_x$-conjugate, then $u = v$.

}
\end{remark}

We demonstrate in Example \ref{2.33} that for $\alpha, \beta \in
\BB_x(R)$ as in (\ref{2.31}), it may happen that $\alpha$ and
$\beta$ are $G$-conjugate and $u \ne v$.

 \begin{example} \label{2.33} {\em
Let $n = p^2$, where $p \ge 3$ is an odd prime. Let $R = \ZZ_n$, let
$\alpha \in \BB_x(R)$ be defined by $\alpha(x) = x + 1$, and let
$\sigma \in G(R)$ be defined by $\sigma(x) = x + px^2$, where we use
elements in $\ZZ$ to represent their equivalence classes in $\ZZ_n$.
Notice that $\sigma^{-1}(x) = x - px^2$. We have:
\begin{eqnarray*}
  \sigma^{-1}\alpha\sigma(x)
    &=&  \sigma^{-1}\alpha(x + px^2)   \\
&=& \sigma^{-1}(x+1  + p(x+1)^2)    \\
&=&  1 + \sigma^{-1}(x) + p\sigma^{-1}(x^2+2x+1)\\
&=& 1 + (x-px^2) + p[(x-px^2)^2 + 2(x-px^2) + 1] \\
&=& 1+p + (1+2p)x.
                      \end{eqnarray*}
Therefore $\beta \in \BB_x(R)$, where $\beta(x) = (1+2p)x + 1+p$ is
$G$-conjugate to $\alpha$ and $1+2p  \ne 1 ~~ \mod p^2$.

}
\end{example}

\begin{remark} \label{2.4} {\em
Let  $\alpha_1$ and $\alpha_2$  be elements  of $G(R)$ defined, as
in (\ref{Gilmer}), by
\begin{eqnarray*}
\alpha_i(x) = c_i + u_i x + x^2 f_i(x),~~~~~i = 1, 2,
\end{eqnarray*}
where $c_i \in R$, $u_i \in \UUU(R)$ and  $f_i(x)  \in \NNN(R[x])$.

\begin{enumerate}
\item   In Example \ref{2.44} we demonstrate  that it is possible
to have $u_1 \ne u_2$ and yet $\alpha_1$ and $\alpha_2$ are
$\BB_x$-conjugate.
\item  In Theorem \ref{2.6} we prove that if
$\alpha_1$ and $\alpha_2$ are in $\BB_x$ and are $G$-conjugate and
if $u_1 = u_2$, then  $\alpha_1$ and $\alpha_2$ are
$\BB_x$-conjugate.
\end{enumerate}

}
\end{remark}

\begin{example} \label{2.44} {\em
Let $n = p^2$, where $p \ge 3$ is an odd prime. Let $R = \ZZ_n$
and let $\alpha \in G(R)$ be defined by $\alpha(x) = x + px^2$,
where we use elements in $\ZZ$ to represent their equivalence
classes in $\ZZ_n$. Let $\sigma \in \BB_x(R)$ be defined by
$\sigma(x) = 1 + x$. Notice that $\sigma^{-1}(x) = -1 + x$. We
have:
\begin{eqnarray*}
  \sigma^{-1}\alpha\sigma(x)
    &=&  \sigma^{-1}\alpha(1 + x)   \\
&=& \sigma^{-1}(1 + x + px^2)    \\
&=&  1 + \sigma^{-1}(x) + p\sigma^{-1}(x^2)\\
&=& 1 + (-1+x) + p(-1+x)^2 \\
&=&  x +p(1 -2x + x^2)\\
&=& p + (1-2p)x + px^2.
                      \end{eqnarray*}
}
\end{example}

\begin{theorem} \label{2.6}
{\it Let  $\alpha$ and $\beta$  be elements  of $\BB_x(R)$ defined by
\begin{eqnarray*}
\alpha(x) = a + u x~,~~~~~~\beta(x) = b + u x,
\end{eqnarray*}
where $u \in \UUU(R)$, and $a, b \in R$. If $\alpha$ and $\beta$ are
$G$-conjugate, then $\alpha$ and $\beta$ are $\BB_x$-conjugate.}
\end{theorem}

\noindent  {\it Proof.} Let $\sigma \in G$ and let $y = \sigma(x)$.
We have
$$
\beta = \sigma^{-1}\alpha \sigma \iff \sigma \beta(x) = \alpha
\sigma(x) = \alpha(y).
$$
Since
$$\sigma\beta(x) = \sigma(b + u x) = b + u\sigma(x) = b + uy,
$$
we have
\begin{eqnarray} \label{2.66}
\alpha(y) = b + uy  \iff \beta = \sigma^{-1}\alpha \sigma.
\end{eqnarray}
The assumption that $\alpha$ and $\beta$ are $G$-conjugate implies
there exists $\sigma \in G$ with $y = \sigma(x)$ and $\alpha(y) = b
+ uy$.  The element   $y$ is of the form
$$
y = d + wx  + f(x) = d + wx + x^2g(x),
$$
where $w \in \UUU(R),~d \in R$ and $f(x)$ and $g(x)$ are nilpotent
element of $R[x]$.

To prove that $\alpha$ and $\beta$ are $\BB_x$-conjugate, it
suffices to find  $Y = c + vx$, where  $v \in \UUU(R)$ and  $c \in
R$, such that  $\alpha (Y) = b + uY$. For, by (\ref{2.66}), if $\tau
\in \BB_x$ is defined by $\tau(x) = Y$, then  $\beta =
\tau^{-1}\alpha \tau \iff \alpha (Y) = b + uY$.

It follows from
\begin{eqnarray*}
b+uy &=& \alpha (y) ~=~ \alpha (wx + d + f(x)) ~=~ w(ux+a) + d + f(ux+a) \\
    &=& u(wx + d + f(x)) + wa - ud - uf(x) + d + f(ux+a)\\
    &=& uy + wa - ud - uf(x) + d + f(ux+a)
                  \end{eqnarray*}
that
\begin{eqnarray*}
   b &=& (1-u)d + wa +  f(ux+a)-u f(x).
                      \end{eqnarray*}
Putting $x=0$ and noting that $f(0)=0$ (since $f(x) = x^2 g(x)$), we
obtain
\begin{eqnarray*}
   b &=& (1-u)d + wa +  f(a) ~=~ (1-u)d + a(w +  a g(a)).
                      \end{eqnarray*}
Therefore
\begin{eqnarray*}
    \alpha (y)
    &=&     u y + b \\
    &=&  uy + (1-u)d + a(w +  a g(a))\\
    &=&  u y + (1-u)d + a v,
                      \end{eqnarray*}
where $v=w + ag(a)$ is a unit since $g(a)$ is nilpotent. Let $Y = vx
+ d$. Then $R[Y] = R[x]$, and
\begin{eqnarray*}
    \alpha (Y) &=& \alpha(vx + d) \\
    &=& v(ux+a) + d \\
    &=& u(vx + d) + va - ud  + d\\
    &=& u Y + (1-u)d + a v\\
    &=& uY +b,
    \end{eqnarray*}
as desired. \hfill $\Box$

\medskip

It is clear that  $\BB_x(R) = G(R)$ if  and only if $R$ is reduced.
One wonders whether there is a more quantitative version of this
statement that relates the relative size of $\NNN (R)$ in $R$ and
that of $\BB_x (R)$ in $G(R)$. In particular, we ask:

\begin{questions} \label{2.2} {\em
Which rings $R$ have the property that
$$ G(R)  =   \bigcup \left\{\BB_z(R) ~~| ~~ z \in R[x] \mbox{ and }
R[x] = R[z] \right\}? $$ If $R$ is not reduced, do there always exist
elements
$$ \xi \in G(R) \setminus \bigcup \left\{\BB_z(R) ~~| ~~ z \in R[x]
\mbox{ and } R[x] = R[z]\right\}? $$ If there exist elements  $ \xi
\in G(R) \setminus \bigcup \left\{\BB_z(R) ~| ~ z \in R[x] \mbox{
and } R[x] = R[z]\right\}$,  what properties distinguish such
elements $\xi$?

}
\end{questions}

\begin{remark} \label{2.21} {\em
If $\sigma \in G(R)$ is defined by $\sigma(x) = z$, then
$\sigma^{-1} \BB_z \sigma = \BB_x$. Therefore every element of
$\BB_z$ is $G$-conjugate to an element of $\BB_x$. Hence if $\xi$ is
in the center of $G(R)$ and $\xi \in \BB_z(R)$ for some $z \in R[x]$
with $R[x] = R[z]$, then
$$\xi \in \bigcap \left\{\BB_z(R) ~ | ~ z \in R[x] \mbox{
and } R[x] = R[z]\right\}.
$$
In Section 3 we observe that for $R = \ZZ_4$, the center of $G(R)$
is not contained in $\BB_x(R)$. We deduce that for $R = \ZZ_4$ we
have
$$ G(R)  \ne   \bigcup \left\{\BB_z(R) ~~| ~~ z \in R[x] \mbox{ and }
R[x] = R[z] \right\}. $$

}
\end{remark}

We examine in detail in Section 3 the structure of the automorphism
group $G(\ZZ_4)$ of $\ZZ_4[x]$ and enumerate the  invariant subrings
of $\ZZ_4[x]$ with respect to subgroups of $G(\ZZ_4)$. We use the
following results that hold in more generality.


\begin{lem} \label{2.71}   Let  $R$  be a ring
and let $f \in R[x]$ be a monic polynomial with $\deg f = d \ge 1$.
For each nonzero polynomial $g(x)  \in R[x]$ there exists an integer
$n \ge 0$ and a unique representation for $g(x)$ as
\begin{eqnarray} \label{decimal}
g(x) = \sum_{k=0}^n   g_k f^k,
\end{eqnarray}
where each  $g_k \in R[x]$, $g_n \ne 0$,   and for each $k$ with $0
\le k \le n$, either $g_k = 0$  or $\deg g_k < d$.

\end{lem}

\noindent  {\it Proof.} If $\deg g < d$, then the statement is clear
with $n = 0$.  We use induction on $\deg g$  and assume for some
integer $m \ge d$ that every polynomial $G$ with $\deg G < m$ can be
represented as in (\ref{decimal}). Let $g \in R[x]$ be a polynomial
with $\deg g = m$ and write $m  = dn + r$, where $n$ and $r$ are
integers and $0 \le r < d$.  Let $c$ denote the leading coefficient
of $g$. Then $G = g - cx^rf^n \in R[x]$ and either $G = 0$
or $\deg G < m$.  If $G = 0$, then $g = cx^rf^n$ has the form given
in (\ref{decimal}). If $G \ne 0$, then by induction $G$,
 and hence also $g = G + cx^rf^n$,  has  the
desired form.

To prove  uniqueness, notice that if $g(x) = \sum_{k=0}^n   g_k
f^k$, where each $g_k$ is either 0 or $\deg g_k < d$,  then $g(x) =
0$ only if all the $g_k$ are 0. For if some $g_k \ne 0$, let $s =
\max\{k ~|~ g_k \ne 0\}$ and let $c$ denote the leading coefficient
of $g_s$. Then $c$ is the leading coefficient of $\sum_{k=0}^n   g_k
f^k$, so this polynomial is nonzero.
  \hfill $\Box$

\begin{theorem} \label{2.8} {\it Let $R$ be a ring and
let $\beta$ be the $R$-automorphism of $R[x]$
defined by $\beta(x) = -x + 1$. Then $\beta^2 = 1$ and the ring of
invariants of the cyclic group $\langle \beta \rangle$ acting on
$R[x]$ is $R[x]^{\langle \beta \rangle} = R[y]$, where $y =x(-x+1)$.}
\end{theorem}

\noindent  {\it Proof.} Every polynomial  $f(x) \in R[x]$ has a
unique representation as
\begin{eqnarray*}
    f(x) &=& \sum_{k=0}^{n} (a_k x + b_k) x^k (x-1)^k,
    \end{eqnarray*}
for some integer $n \ge 0$,  where $a_k$ and $b_k$ are in $R$, $1
\le k \le n$. Assume that $f(x)$ is fixed by $\beta$.  Then $f(x) =
f(1-x)$,  and therefore
\begin{eqnarray*}
    \sum_{k=0}^{n} (a_k x + b_k) x^k (x-1)^k  &=&
    \sum_{k=0}^{n} (-a_k x + a_k + b_k) x^k (x-1)^k.
    \end{eqnarray*}
By uniqueness of representation, we conclude that $a_k  = 0$ for
every $k$.  \hfill $\Box$

\begin{remark} \label{2.99} {\em
Let $R$ be a ring and
let
$$
f = a_0 + a_1x + \cdots + a_nx^n \in R[x]
$$
be a polynomial. Consider the following assertions:
\begin{enumerate}
\item The surjective $R$-algebra homomorphism of $R[x]$ onto $R[f]$
defined by mapping $x \mapsto f$ is injective.
\item The subring $R[f]$ of $R[x]$ is $R$-isomorphic to $R[x]$.
\item The $R$-algebra $R[f]$ is a polynomial ring over $R$. \item
The annihilator in $R$  of the ideal $I = (a_1, \ldots, a_n)R$ is
zero.
\end{enumerate}
It is readily seen that the first 3 assertions  are equivalent, and
it is well known that these are also equivalent to assertion 4
\cite[Theorem 2]{Gilmer}. In considering condition 4, since $R[f] =
R[f - a_0]$,  one may assume that $a_0 = 0$.  With this assumption,
if $b \in R$ is nonzero and $bI = 0$, then $bx \in R[x]$ is a
nonzero polynomial in the kernel of the $R$-algebra homomorphism
defined by $x \mapsto f$. On the other hand, if this $R$-algebra
homomorphism is not injective, let
\begin{eqnarray*}
b_mx^m + \cdots + b_1x + b_0 \in R[x]
\end{eqnarray*}
be a nonzero polynomial of minimal degree  such that
\begin{eqnarray*}
b_mf^m + \cdots + b_1f + b_0  = 0.
\end{eqnarray*}
Then $a_0 = 0$ implies $b_0 = 0$. Therefore
\begin{eqnarray*}
(b_mf^{m-1} + \cdots + b_1)f = 0,
\end{eqnarray*}
and $b_mf^{m-1} + \cdots +  b_1$ is a nonzero polynomial  because of
the minimal degree assumption. Thus
$$f = a_nx^n + \cdots + a_1x
$$ is
a zero-divisor in $R[x]$. A well known theorem of McCoy \cite[page
290]{Mc} implies that the ideal $I$ of $R$ has a nonzero annihilator
(see also \cite[Theorem~4, pages 34-36]{Mc2},   \cite[page
330]{Gilmer} and \cite[Exercise 2(iii), page 11]{AM}).

 }
\end{remark}

\begin{theorem} \label{2.9}
 Let $R$ be a ring and let  $f \in R[x]$ be such that $R[f]$ is a
polynomial ring over $R$. If $g$ is a  nilpotent element of
$R[x]$, then the ring $R[f,g]$ is a polynomial ring over $R$ if
and only if $g \in R[f]$.
\end{theorem}

\noindent {\it Proof.} Let $N=\NNN (R)$ be the nilradical of $R$.
The nilradical of the polynomial ring $R[x]$ is given by $\NNN
(R[x]) = N R[x]$. Let $R/N = F$. Then $R[x]/(NR[x]) = F[x]$.

We denote the image of $r \in R$ under the canonical map $R
\rightarrow R/N$ by $\overline{r} = r + N$. This map extends to a
map from $R[x]$ to $F[x]$.  We denote the image of $f \in R[x]$
under  this map by $\overline{f} = f + N[x]$. Thus if $f = \sum
r_i x^i \in R[x]$, then $\overline{f} = \sum \overline{r_i} x^i$.

Assume that $R[f,g] = R[y]$ for some $y$ in $R[x]$. Passing to
quotients mod $N R[x]$, we see that
$$F[\overline{y}] = F[\overline{f},\overline{g}] = F[\overline{f}].$$
By Fact \ref{1.00}, we have
$$\overline{y} = U \overline{f} + C,$$
where $U$ is a unit in $F$ and $C \in F$. Thus $U = u + N$, $C = c + N$,
and $(u+N)(v+N)=1+N$ for some $v \in R$. Hence $uv=1+\nu$ for some $\nu \in N$
and $uv$ and therefore $u$ is a unit in $R$. Hence
$$y = u f + c + H(x),$$
where $H(x) \in NR[x]$. Again, by Fact \ref{1.00}, $R[y] = R[f]$,
as desired. \hfill $\Box$

\section{Automorphisms of $\ZZ_4[x]$}

In this section, we determine the structure of the group $G(R)$ of
$R$-automorphisms of the polynomial ring $R[x]$ in the case where
$R=\ZZ_4$ and we describe the ring of invariants $R[x]^H$ for
every subgroup $H$ of $G(R)$. We also describe the conjugacy
classes in $G(R)$.

\begin{remark} \label{4.0} {\em
Let $R$ be a ring and let $f \in \NNN(R[x])$ be a nilpotent
element of the polynomial ring $R[x]$. We associate with $f$ the
$R$-automorphisms $\alpha_f$ and $\beta_f$ in $G(R)$ defined by
\begin{eqnarray}
    \alpha_f : x \mapsto  x +  f~,~~~~~~
\beta_f : x \mapsto  - x + 1 +  f. \label{alphafbetaf}
\end{eqnarray}
Notice that the correspondences $f \mapsto \alpha_f$ and $f
\mapsto \beta_f$ are both one-to-one. Moreover, the set $\AAA :=
\{\alpha_f : f \in \NNN(R[x])\}$ is a subgroup of $G(R)$. Also the
constant term of $\alpha_f(x)$ is a nilpotent element of $R$,
while the constant term of $\beta_f(x)$ is a unit of $R$ for each
$f \in \NNN(R[x])$. Therefore the sets $\AAA$ and $\{\beta_f : f
\in \NNN(R[x])\}$ are disjoint. In the special case where $2R$ is
a maximal ideal of $R$ with $R/2R = \ZZ_2$, so, in particular, in
the case where $R = \ZZ_4$, every automorphism\footnote{In the
case where $R = \ZZ_4$, every automorphism of $R[x]$ is an
$R$-automorphism.} of $R[x]$ is of the form $\alpha_f$ or
$\beta_f$ for some $f \in \NNN(R[x])$. Indeed, in this case,
$G(R/2R)$ is a group of order 2, and $\AAA$ is the normal subgroup
of $G(R)$ that is the kernel of the canonical surjective
homomorphism of $G(R)$ onto $G(R/2R)$ while $\{\beta_f : f \in
\NNN(R[x])\}$ is the unique nonidentity coset of $\AAA$ in $G(R)$.

 }
\end{remark}

We start by describing the ring of invariants $R[x]^H$ in the case
where $R = \ZZ_4$ and  $H$ is a cyclic subgroup of $G(R)$. By
Remark \ref{4.0}, $H$ has the form $H = \langle \alpha_f\rangle$
or $H = \langle \beta_f \rangle$ for some $f \in \NNN(R[x])$. The
fixed rings of each of these types is described in Theorems
\ref{4.1}  and \ref{4.3} below. Lemma \ref{4.2} is used in the
proof of Theorem \ref{4.3}.

\begin{theorem} \label{4.1}   Let $R = \ZZ_4$ and let $\alpha \in G(R)$ be defined by
$$\alpha (x) = x + f,$$
where $f \in \NNN(R[x])$.  If $f = 0$, then the order of $\alpha$
is 1. If $f \ne 0$, then the order of $\alpha$ is  2 and the fixed
ring of $\alpha$ is $R[x^2,2x]$ and is not a polynomial ring.
\end{theorem}

\noindent {\it Proof.} It is clear that $\alpha$ is the identity
element of $G(R)$ if and only if $f = 0$. Assume that $f \ne 0$.
Since  $\NNN(R[x]) = 2R[x]$, we have $2f = f^2 = 0$. It follows
that $\alpha^2(x) = x$ and $\alpha$ has order 2. Also $f = 2g$,
where $g$ has the form, for some integer $m \ge 0$,
$$g = \sum_{k=0}^m b_k x^k,~~ \mbox{ where } ~~ 2b_m \ne 0.$$
Let $S = R[x^2,2x]$.  Clearly, $S \subseteq R[x]^{\langle \alpha
\rangle}$. To show that this inclusion is an equality,  assume
that there exists an element $h \in R[x]^{\langle \alpha \rangle}
\setminus S$.  Subtracting from $h$ an element in $S$,  we obtain
for some positive integer $n$ an element $h'  \in R[x]^{\langle
\alpha \rangle} \setminus S$ of the form
$$
\sum_{k=0}^n a_k x^{2k+1}, ~~\mbox{ where } a_n = 1.
$$
 Since $h'$
is $\alpha$-fixed, it follows that
\begin{eqnarray*}
 \sum_{k=0}^n a_k x^{2k+1}    &=& \sum_{k=0}^n a_k x^{2k}(x+2g) \\
\sum_{k=0}^{n-1} a_k x^{2k+1}    &=& 2 x^{2n}g
+\sum_{k=0}^{n-1}a_k x^{2k}(x+2g).
    \end{eqnarray*}
Comparing the coefficients of the highest degree terms in this last
equation, we see that $2b_mx^{2n+m} = 0$.  This contradicts the
assumption that $2b_m \ne 0$. Therefore $R[x]^{\langle \alpha
\rangle} = S$. By Theorem \ref{2.9}, $S = R[x^2, 2x]$ is not a
polynomial ring over $R$\hfill $\Box$

\begin{lem} \label{4.2}
 Let $R = \ZZ_4$ and let $\theta$ be the automorphism of $R[x]$
defined by
$$\theta (x) = x + 1.$$
Then $$R[x]^{\langle \theta \rangle} = R[w^2,2w],$$ where
$w=x(x+1)$. In particular,  $R[x]^{\langle \theta \rangle}$ is not
a polynomial ring over $R$.
\end{lem}

\noindent {\it Proof.} Let $M =  x(x+1)(x+2)(x+3)$ be the norm of
$x$ with respect to $\theta$.  Then $\theta (w) = (x+1)(x+2) = w +
2(x+1).$ Therefore, both $w^2$ and $2w$ are fixed by $\theta$. It
remains to show that every $g \in R[x]^{\langle \theta \rangle}$
belongs to $R[w^2,2w]$. By Lemma \ref{2.71}, the polynomial  $g$
has a unique representation as
\begin{eqnarray*}
g(x) = \sum_{k=0}^n   g_k M^k,
\end{eqnarray*}
for some integer $n \ge 0$,  where each  $g_k$ is either 0 or a
polynomial in $R[x]$ of degree  less than 4. Since $\theta$ fixes
$M$ and does not increase degrees, it follows from Lemma
\ref{2.71} that $\theta (g) = g$ if and only if $\theta
(g_k) = g_k$ for all $k$. Thus it suffices to show that if
 $h = ax+bx^2+cx^3 \in R[x]$ is fixed by $\theta$, then $h \in R[w^2,2w]$.
 We have
\begin{eqnarray*}
    \theta (h) = h
    &\ifff&
ax+bx^2+cx^3\\
&& =
a(x+1)+b(x^2+2x+1)+c(x^3+3x^2+3x+1)\\
&\ifff&
3c = 3c + 2b = a+b+c = 0\\
&\ifff&
c = 0, ~a = -b, ~2b = 0\\
&\ifff&
h = b(x^2-x), ~2b = 0\\
&\ifff&
h = 2d (x^2-x) = 2d(x^2+x) , \mbox{~for some $d \in R$}\\
&\ifff& h = 2dw.
\end{eqnarray*}
Therefore $R[x]^{\langle \theta \rangle} = R[M,2w] =R[w^2,2w]$,
since $M=w^2+2w.$ By Theorem \ref{2.9}, $R[x]^{\langle \theta
\rangle}$ is not a polynomial ring over $R$. \hfill $\Box$

\begin{theorem} \label{4.3}  {\it Let $R = \ZZ_4$ and let $\beta \in G(R)$ be
defined by
$$\beta (x) = - x + 1 + f,$$
where $f \in \NNN(R[x])$.  Let $y = x(-x+1)$.   The order
of $\beta$ is either 2 or 4, and the following are equivalent:
\begin{enumerate}
\item  The order of $\beta$ is 2. \item The element $f \in R[y]$.
\item $f = 2h$ for some $h \in R[y]$. \end{enumerate}
If the order
of $\beta$ is 2, then the fixed ring of $\beta$ is $R[y+xf]$, a
polynomial ring over $R$ generated by the element $y+xf$. If the
order of $\beta$ is 4, then the fixed ring of $\beta$ is
$R[y^2,2y]$ and is not a polynomial ring.}
\end{theorem}

\noindent {\it Proof.} It is clear that
\begin{eqnarray}
    \beta^2 (x) = x + g, \label{sigma2}
\end{eqnarray}
where $g =  \beta (f) -f$. Since $\beta(\NNN(R[x]) = \NNN(R[x])$,
$\beta(f) - f  \in \NNN(R[x])$. By Theorem \ref{4.1}, the order of
$\beta^2$ is 1 or 2. Since $\beta$ is not the identity element of
$G(R)$, the order of $\beta$ is 2 or 4.

We consider first the case where the order of $\beta$ is 2.
Clearly,
\begin{eqnarray}
\mbox{order ($\beta$) = 2~}&\ifff& \beta^2(x) = x ~\ifff~ \beta
(f)- f = 0. \label{order2}
\end{eqnarray}
We show  that this happens if and only if $f = 2h$ for some $h \in
R[y]$.

It is easy to see that $\beta (x^k) = (-x+1)^k$ if $k$ is even and
$\beta (2x^k) = 2(-x+1)^k$ if $k$ is odd. Thus letting $\beta_0$
be the automorphism $x \mapsto -x+1$, we see that $\beta$ and
$\beta_0$ coincide on $2x^k$ for every $k$ and therefore coincide
on every element of $\NNN(R[x])$. Therefore
\begin{eqnarray*}
\beta (f) -f  = \beta_0 (f) -f.
\end{eqnarray*}
 From this and (\ref{order2}) it follows that the order of $\beta$ is
2 if and only if $f$ belongs to $R[x]^{\langle \beta_0 \rangle}$. By
Theorem \ref{2.8},  $R[x]^{\langle \beta_0 \rangle} = R[y]$.
Therefore
\begin{eqnarray*}
\mbox{order ($\beta$) = 2~}&\ifff& f \in R[y].
\end{eqnarray*}
However an element in $R[y]$ that  2 multiplies to 0 must be of the
form $2h$ for some $h \in R[y]$. Therefore
\begin{eqnarray*}
\mbox{order ($\beta$) = 2~}&\ifff& f = 2h \mbox{~for some~} h \in
R[y].
\end{eqnarray*}

Assume  that the order of $\beta$ is 2. Thus $f = 2h \in 2R[y]$.
Notice that $\beta (y) = y + 2h$, and therefore $\beta (y^2) =
y^2$ and $\beta (2y) = 2y$. Therefore $\beta (2y^k) = 2y^k$ for
all $k$ and hence $\beta (2g) = 2g$ for all $g \in R[y]$. In
particular, this holds for $f =2h$, and we have
\begin{eqnarray*}
    \beta (2xh) &=& (-x+1+2h)(2h)\\
    &=& 2xh + 2h.
\end{eqnarray*}
Let $z = x + 2xh.$ By Fact \ref{1.00}, $R[x] = R[z]$. Also,
\begin{eqnarray*}
    \beta (z) &=& (-x+1+2h)+(2xh+2h)\\
    &=& -z + 1.
\end{eqnarray*}
Therefore
\begin{eqnarray*}
R[x]^{\langle \beta\rangle} &=& R[z]^{\langle \beta\rangle}\\
&=&  R[z(-z+1)] \mbox{~by Lemma \ref{2.8}}\\
&=& R[(x+2xh)(-x-2xh+1)] \\
&=& R[x(-x+1) + 2xh]\\
&=& R[y + 2xh]\\
&=& R[y +xf],
\end{eqnarray*}
as claimed. Note that $R[y + xf]$  contains $R[y^2,2y]$ properly
since $R[y^2,2y]$  is not a polynomial ring over $R$ by Theorem
\ref{2.9}.

Assume that the order of $\beta$ is 4. Then the order of $\beta
^2$ is 2. Since $\beta^2 (x)$ has the form given in
(\ref{sigma2}), Theorem \ref{4.1} implies that  $R[x]^{\langle
\beta^2 \rangle} = R[x^2,2x]$. This does not depend on $f$. Nor
does the action of $\beta$ on this ring, since
$$\beta (x^2) = (-x+1)^2~,~~\beta (2x) = 2x + 2.$$
Thus for the purpose of finding $R[x]^{\langle \beta \rangle}$,
one may take $f =2x$. Then  $\beta : x \mapsto x+1$. By Lemma
\ref{4.2}, $R[x]^{\langle \beta \rangle} = R[w^2,2w]$, where $w =
x(x+1)$, Also $R[w^2, 2w] = R[y^2, 2y]$, so  $R[x]^{\langle \beta
\rangle} = R[y^2,2y]$.  This completes the proof of Theorem
\ref{4.3}. \hfill $\Box$

\begin{notation} \label{4.4} {\em
To  consider $R[x]^H$ for an arbitrary subgroup $H$ of
$G(R)$, we use the following notation. Referring to
(\ref{alphafbetaf}), we observe that $\alpha_0$ is the identity of
$G(R)$ and $\beta_0$ is the automorphism $\beta$ defined by
\begin{eqnarray*}
\beta := \beta_0 : x \mapsto -x+1,
\end{eqnarray*}
that is  considered  in Theorem \ref{2.8}. The automorphism
$\theta: x \mapsto x+1 $ of  Lemma  \ref{4.2} is $\beta_{2x}$ and
its inverse is $\beta_{2(-x+1)}$. We also let $y$ be defined by
\begin{eqnarray}
    y = x(-x+1),
\end{eqnarray}
as in Theorem \ref{4.3}. For $h \in R[x]$, we denote $\beta_0 (h)$
by $h'$. Thus $h'$ is obtained from $h$ by replacing $x$ by
$-x+1$. Hence $h'' = h$ for all $h \in R[x]$ and
\begin{eqnarray}
    h = h' &\ifff& h \in R[x]^{\langle \beta \rangle} = R[y].
\end{eqnarray}
Also the map  $\phi : R[x] \rightarrow R[y]$ defined by $\phi (f) = f+f'$ is onto. In fact,
\begin{eqnarray*}
\phi \left(
x^{n+1} (-x+1)^n\right) =
x^{n+1} (-x+1)^n+x^{n} (-x+1)^{n+1} &=& x^n (-x+1)^n \\
&=& y^n.
\end{eqnarray*}

}
\end{notation}

Theorem \ref{4.5}  describes $R[x]^H$ for an arbitrary subgroup
$H$ of $G(R)$. These fixed subrings are precisely those subrings
fixed by cyclic subgroups (as described in Theorems \ref{4.1} and
\ref{4.3}), together with the family of polynomial rings
$R[y+xf],~f\in R[y]$, where $y = x(-x+1)$. We let $e$ denote the
identity element of the group $G(R)$.

\begin{theorem} \label{4.5}  Let $R = \ZZ_4$ and let $\AAA =
\{\alpha_f : f \in \NNN(R[x])\}$ be the normal subgroup of $G(R)$
of index 2 defined  in Remark \ref{4.0}. Let $H$ be a subgroup of
$G(R)$, and let  $y = x(-x+1)$.
\begin{enumerate}
\item[(a)]  If $H$ is a subgroup of $\AAA$, then $R[x]^H$ is
either $R[x]$ or  $R[x^2,2x]$ depending on whether or not $H$ is
trivial. \item[(b)] If $H$ is not contained in  $\AAA$ and if $H
\cap \AAA \ne \langle e \rangle$,
        then $R[x]^H = R[y^2,2y]$.
\item[(c)]    If $H$ is not contained in  $\AAA$ and if $H \cap
\AAA = \langle e \rangle$,  then $H$ is cyclic generated by an
element $\beta_f$ of order 2 and $R[x]^H = R[y + xf]$, where $f
\in \NNN(R[y])$.
\end{enumerate}

\end{theorem}

\noindent {\it Proof.} (a) Theorem \ref{4.1} implies that if $H$
is a non-trivial subgroup of $\AAA$, then $R^H  = R[x^2,2x]$.

(b) Assume that $H$ is not contained in  $\AAA$ and that $H_0 : =
H \cap \AAA$ is nontrivial. By (a), $R[x]^{H_0} = R[x^2,2x].$
     Therefore $R[x]^H \subseteq  R[x^2,2x]$. Also,
     every $\beta_h \in H$ acts on $R[x^2,2x]$ as follows:
     $$x^2 \mapsto (-x+1)^2~,~~2x \mapsto 2(-x+1).$$
     Therefore
every $\beta_h \in H$ fixes $y^2$ and $2y$ and hence
\begin{eqnarray}
    R[x]^H \supseteq R[y^2,2y]. \label{supset}
\end{eqnarray}
Letting $\alpha$,  $\beta$, and $\theta$ be the automorphisms on $R[x]$ defined by
$$\alpha : x \mapsto -x~,~~\beta :
x \mapsto -x +1~,~~\theta = \alpha \beta : x \mapsto x+1,$$
 we see that $\beta$ restricts to an automorphism of $R[x^2,2x]$,  and
\begin{eqnarray*}
R[x]^H &=& R[x^2,2x]^{\langle \beta \rangle}\\
&=&
\left(R[x]^{\langle \alpha \rangle}\right)^{\langle \beta \rangle}\\
&\subseteq&
R[x]^{\langle \alpha \beta \rangle}\\
&=&
R[x]^{\langle \theta \rangle}\\
&=&
R[y^2,2y] \mbox{~(by Theorem \ref{4.2})}.
\end{eqnarray*}
 From this and (\ref{supset}) it follows that $R[x]^H = R[y^2,2y]$, as claimed.

    (c) Assume that $H$ is not contained in  $\AAA$ and that $H \cap \AAA =
    \langle e\rangle $.
    Since $[G;\AAA] = 2$, we have $[H:H \cap \AAA] \le 2$.
    Therefore $H$ is a cyclic group of order 2, and $H = \langle \beta_g \rangle$,
    where $\beta_g^2 = 1$.  By Theorem \ref{4.3}, $g \in \NNN(R[y])$ and
    $R[x]^H = R[y+xg]$.
            This completes the proof of Theorem \ref{4.5}.  \hfill $\Box$

\noindent \begin{remark} \label{4.6}  {\em
Theorem \ref{4.5}
asserts that every  ring of invariants of  $R[x]$ with respect to
a subgroup of $G(R)$ is one of the following rings:
    $$\mbox{(i)}~ R[x], ~~~~~\mbox{(ii)}~ R[x^2,2x], ~~~~~\mbox{(iii)}~  R[y^2,2y]~
    ,~~~\mbox{ or }~~~~~\mbox{(iv)} ~ R[y+xg],  $$
    for some  $g \in \NNN(R[y])$. The first three items are
    specific rings while item (iv) describes an infinite family of
    polynomial subrings of $R[x]$. The ring $R[y^2, 2y]$ is the
    ring of invariants of $R[x]$ with respect to $G(R)$ and thus
    is the unique smallest ring in the family. The rings that are
    fixed rings with respect to subgroups of $\BB_x(G)$  are the
    rings of the first three items and the polynomial rings
    $R[y]$ and $R[x(x+1)]$.  Notice that $R[y]$ corresponds to $g = 0$ in
    (iv),
    and $R[x(x+1)]$  corresponds to $g = 2$ or $g = 2x$. Letting
    $$\alpha_{2x} : x \mapsto -x ~, ~~
    \beta_0 : x \mapsto -x+1~,~~\beta_2 : x \mapsto -x-1 ~, ~~ \theta : x \mapsto x+1~,$$
we easily see that
    \begin{eqnarray*}
    R[x] &=& R[x]^{\langle e\rangle} \\
    R[x^2,2x] &=& R[x]^{\langle \alpha_{2x} \rangle}\\
    R[y] &=& R[x]^{\langle \beta_0 \rangle}\\
    R[x(x+1)] &=& R[x]^{\langle \beta_2 \rangle}\\
    R[y^2,2y] &=& R[x]^{\langle \theta \rangle}
\end{eqnarray*}

             }
\end{remark}

   The polynomial subrings    $R[y+xg]$
      of item (iv) for $g \in \NNN(R[y]) \setminus \{0,2, 2x \}$
are not rings of
    invariants of  subgroups of $\BB_x$.
However, each of these rings is the fixed ring of an element $\sigma
\in \BB_z(R)$ for some $z$ such that $R[z] = R[x]$. To see this,
simply take $z = x + x g$ and $\sigma : z \mapsto -z + 1$. Then
$R[z] = R[x]$ because $g$ (and hence $xg$) is nilpotent,
$$R[x]^{\langle \sigma \rangle} = R[z]^{\langle \sigma \rangle} =
R[z(-z+1)],$$ and
$$z(-z+1) = (x + xg)(-x+1+xg) = x(-x+1)+xg = y+xg,$$ as desired.

    The rings $R[y+xg]$, $g \in \NNN (R[y]$, are also pairwise different. To show this, we
 associate  with each  subring $S$ of
    $R[x]$, the subgroup
        $$S^* = \{\sigma \in G(R)~ |~ \sigma(s) = s \mbox{ for all }
        s
    \in S \}
    $$
    of $G(R)$     consisting of the  automorphisms that restrict to
    the identity map on $S$. For subrings $S_1$ and $S_2$ of
    $R[x]$, it is clear that  $S_1^* \ne S_2^*$  implies $S_1 \ne S_2$.
    Therefore Theorem \ref{4.8} implies that the rings of
    invariants enumerated in Remark \ref{4.6} are all distinct.

     \begin{theorem} \label{4.8}
      For a subring $S$ of $R[x]$, let $S^*$ denote the
subgroup of $G(R)$ whose elements fix $S$. Then
    \begin{eqnarray*}
    R[x]^* &=& \{ e\}\\
    R[x^2,2x]^* &=& \AAA\\
       R[y^2,2y]^* &=& G(R)\\
    R[y+xf]^* &=& \langle \beta_f \rangle.
    \end{eqnarray*}
    for each  $f \in \NNN(R[y])$.

     \end{theorem}

\noindent {\it Proof.} The first equality is clear. The second
equality follows because $x^2$ and $2x$ are fixed by $\alpha_f$ for
all $f$ and  no $\beta_f$ fixes $2x$ since
\begin{eqnarray*}
    \beta_f (2x) = -2x + 2 = 2x + 2 \ne 2x.
\end{eqnarray*}
Since  $y^2$, and $2y$ are fixed by $\alpha_f$ and $\beta_f$ for
all $f$, the third equality is clear. To establish the last
equality, observe that
\begin{eqnarray*}
    \alpha_g (y + xf) = y+xf &\ifff& \alpha_g (x(-x+1) + xf) = x(-x+1)+xf\\
     &\ifff& (x+g)(-x-g+1) + (x+g)f = x(-x+1)+xf\\
&\ifff& x(-x+1)+g+xf = x(-x+1)+xf\\
    &\ifff& g = 0.\\
        \beta_g (y+xf) = y + xf &\ifff& \beta_g (x(-x+1) + xf) = x(-x+1)+xf\\
                &\ifff& (-x+1+g)(x+g) + (-x+1+g)f    = x(-x+1)+xf \\
                &\ifff& x(-x+1) + g  - xf + f + gf  = x(-x+1)+xf\\
                &\ifff&  g   + f  = 0\\
&\ifff&  g = f.
\end{eqnarray*}
This completes the proof of Theorem \ref{4.8}.   \hfill $\Box$

\medskip

We now turn to the structure and conjugacy classes of $G(R)$, again for $R=\ZZ_4$.

\begin{theorem} \label{4.9}
{\it Let $R = \ZZ_4$ and let $\alpha_f$ and $\beta_f$ be as defined
in (\ref{alphafbetaf}). Let $y=x(-x+1)$,  and let
\begin{eqnarray*}
    \AAA = \{ \alpha_f : f \in \NNN(R[x])\},~
    \AAA_0 = \{ \alpha_f : f \in \NNN(R[y])\},~
            \CCC = \{ \alpha_0 = e, \beta_0 = \beta \}.
    \end{eqnarray*}
Then the groups $\AAA$ and $\CCC$ are isomorphic to the additive groups $\ZZ_2[x]$
and $\ZZ_2$, respectively,  and  $G(R)$ is the
extension of $\AAA$ by $\CCC$ via the multiplication
$$\beta^{-1} \alpha_g \beta = \beta_{g'},$$
where $g' = \beta_0(g)$. Also the center of $G(R)$ is $\AAA_0$.
}
 \end{theorem}

\noindent
{\it Proof.}
It is clear that   $\alpha_f$ and $\beta_f$ act on $2 R[x]$ and
that
\begin{eqnarray}
    \alpha_f (2h) = 2h~,~~
    \beta_f (2h) = \beta_0 (2h) = 2h'
    \end{eqnarray}
for all $f, h \in R[x]$. Thus the actions of $\alpha_f$ and $\beta_f$ on
$2 R[x]$ are independent of  $f$. It is also easy to verify that
\begin{eqnarray} \label{multiplication}
    \alpha_g \alpha_h = \alpha_{g+h}~,~~
     \beta_g \beta_h = \alpha_{g+h'}~,~~
     \alpha_g \beta_h =   \beta_{g+h}    ~,~~
     \beta_g \alpha_h = \beta_{g+h'},
     \end{eqnarray}
and therefore
\begin{eqnarray} \label{conjugation}
    \alpha_g^{-1} \alpha_h  \alpha_g = \alpha_h~,~~
     \beta_g^{-1} \alpha_h  \beta_g = \alpha_{h'}~,~~
        \alpha_g^{-1} \beta_h  \alpha_g = \beta_{g'+h+g}~,~~
        \nonumber \bigskip  \\
    \beta_g^{-1} \beta_h  \beta_g = \beta_{g+h'+g'}. \hskip 1.4 in
     \end{eqnarray}
The assertions in Theorem \ref{4.9} follow immediately from (\ref{multiplication}) and
(\ref{conjugation}). \hfill $\Box$

\begin{remark} \label{4.10} {\em
To describe conjugacy classes in $G(R)$, let
\begin{eqnarray*}
    \AAA &=& \{ \alpha_f : f \in \NNN(R[x])\},~
    \AAA_0 ~=~ \{ \alpha_f : f \in \NNN(R[y])\}\\
        \BBB &=& \{ \beta_f : f \in \NNN(R[x])\},~~
    \BBB_0 ~=~ \{ \beta_f : f \in \NNN(R[y])\}.
\end{eqnarray*}
We see that $G$ is the disjoint union of $\AAA$ and $\BBB$, that
$\AAA$ is a normal subgroup of $G$ of index 2, and that $\AAA_0$ is
the center of $G$. We also see that each of $\AAA_0$, $\AAA
\setminus \AAA_0$, $\BBB_0$, $\BBB \setminus \BBB_0$ is the union of
conjugacy classes. In $\AAA_0$, a conjugacy class has one element.
In $\AAA \setminus \AAA_0$, a conjugacy class has two elements
$\alpha_f$ and $\alpha_{f'}$. In $\BBB$, two elements $\beta_g$ and
$\beta_h$ are  conjugate if and only if $g -h$ or $g-h'$ belongs to
$R[y]$. In particular, $\BBB_0$ is a conjugacy class. Identifying
$\BBB$ (as a set) with $R[x]$, and $\BBB_0$ with $R[y]$, the
conjugacy class $\BBB_0$ corresponds to $R[y]$ and every other
conjugacy class in $\BBB$ corresponds to the union of two cosets in
the group $R[x]/R[y]$.

Since the center $\AAA_0$ of $G(R)$ is not contained in $\BB_x(R)$,
Remark \ref{2.21} implies that
$$ G(R)  \ne   \bigcup \left\{\BB_z(R) ~~| ~~ z \in R[x] \mbox{ and }
R[x] = R[z] \right\}. $$

}
\end{remark}

\begin{remark} \label{4.11} {\em
In Theorem \ref{2.6} we prove for a general ring $R$ that certain
elements of $\BB_x(R)$ that are conjugate as elements of  $G(R)$ are
actually conjugate as elements of $\BB_x(R)$. For $R = \ZZ_4$, we
show this holds without any additional condition on the elements;
that is,  if two elements in $\BB_x(R)$ are $G(R)$-conjugate, then
they are $\BB_x(R)$-conjugate. This is equivalent to showing that
\begin{eqnarray} \label{4.13}
&&\mbox{$\forall ~ \sigma \in \BB_x(R),~~~~[\sigma]_{\BB_x} =
[\sigma]_{G} \cap \BB_x$,}
\end{eqnarray}
 where
$[\sigma]_{\BB_x}$ and $[\sigma]_{G}$ are the conjugacy classes in
$\BB_x(R)$ and in $G(R)$ that contain $\sigma$.

}
\end{remark}

The 8 elements of $\BB_x(R)$ are:
$$\alpha_0 = e : x \mapsto x~,~~ \alpha_2 : x \mapsto x+2~,~~
\alpha_{2x} : x \mapsto -x~,~~
\alpha_{2x+2} : x \mapsto -x+2,$$
$$\beta_0 : x \mapsto -x +1~,~~
\beta_2 : x \mapsto -x -1~,~~
\beta_{2x} : x \mapsto x +1~,~~
\beta_{2x+2} : x \mapsto x -1.$$

We  verify Condition \ref{4.13}  by observing that each of the pairs
$$(\alpha_{2x}, \alpha_{2x+2})~,~~
(\beta_0, \beta_2)~,~~
(\beta_{2x}, \beta_{2(x+2)})$$
consists of $\BB_x$-conjugate elements. It is direct to check that
$$
\beta_{2x}^{-1}\alpha_{2x}\beta_{2x}=\alpha_{2x+2}~,~~
\alpha_{2x}^{-1}\beta_{0}\alpha_{2x}=\beta_{2}~,~~
\alpha_{2x}^{-1}\beta_{2x}\alpha_{2x}=\beta_{2x+2}.$$

The question of which elements  $\xi \in G(R)$
do  not belong to any $\BB_z(R)$ with $R[z] = R[x]$  is answered
for $R = \ZZ_4$ as follows:
These are precisely those $\xi$ for which $[\xi]_G$ does not
intersect $\BB_x$. These are

(i)  All $\alpha_f$ that do not belong to $\BB_x$, i.e., all
$\alpha_f$ except
$$\alpha_0 =e~,~~ \alpha_2 : x \mapsto x+2~,~\alpha_{2x} :
x \mapsto -x~,~~\alpha_{2x+2} : x \mapsto -x+2.$$

(ii) All $\beta_f$ except
$$\beta_0 : x \mapsto -x+1~,~~ \beta_2 : x \mapsto -x-1~,~~\beta_g~,~~
\beta_{2x + g}~, ~~\beta_{2(x+2)+g},$$
where $g \in \NNN (R[y]$.

\section{Conjugacy classes in the group $\BB(\ZZ_n)$}

Throughout this section, we fix a natural number $n > 1$. All
numbers are elements in $\ZZ$ and an element in $\ZZ_n$ is
represented by one of its inverse images under the natural map $\ZZ
\rightarrow \ZZ_n$. In particular, if $\alpha$ is an
automorphism\footnote{Recall that every automorphism of $\ZZ_n[x]$
is a $\ZZ_n$-automorphism.}  of the polynomial ring $\ZZ_n[x]$ that
is $x$-basic, then $\alpha (x)$ has the form $ux + a$, where $u, a
\in \ZZ$ and where $u$  is a unit mod $n$. We denote the group of
$x$-basic $\ZZ_n$-automorphisms of $\ZZ_n[x]$ by $\BB(\ZZ_n)$.

If $a, b \in \ZZ$,  then $(a,b)$ denotes  the positive greatest
common divisor of $a$ and $b$. Assume  the elements $\alpha$ and
$\beta$ of $\BB(\ZZ_n)$ are such that
\begin{eqnarray}  \label{ab}
    \alpha (x) = ux + a, ~~~~~~~~\beta (x) = vx + b,
\end{eqnarray}
where $(u,n) = (v,n) = 1$. We emphasize the trivial fact that
\begin{eqnarray*}
    \alpha = \beta &\ifff& ux + a = vx + b \mbox{~in~} \ZZ_n[x]\\
    &\ifff& u \equiv v ~(\mod n)  \mbox{~and~}  a \equiv b ~(\mod n).
    \end{eqnarray*}
On the other hand, $\alpha$ and $\beta$ are said to be {\it
equivalent}, and we write $\alpha \cong \beta$,  if they are
conjugate as elements in $\BB(\ZZ_n)$. This happens if and only if
there exists $X = wx + c$, where  $(w,n) =1$,  such that $\alpha (X)
= v X + b$. We have
\begin{eqnarray*}
\alpha (X) = vX + b &\ifff& w(ux+a) + c  = v(wx+c) +b  \mbox{~in~} \ZZ_n[x]\\
    &\ifff& wux+ wa + c  = vwx+ vc +b   \mbox{~in~} \ZZ_n[x]     \\
    &\ifff& wu \equiv vw  ~(\mod n)  \mbox{~~and~~}  wa+c  \equiv  vc +b ~(\mod n) \\
    &\ifff& u \equiv v  ~(\mod n)  \mbox{~~and~~}  wa \equiv (v-1)c + b ~(\mod n).
\end{eqnarray*}
We record the conclusion as:

\begin{fact} \label{3.0} If $\alpha$ and $\beta$ are defined as in (\ref{ab}),
then
\begin{eqnarray*}
    \alpha \cong \beta &\ifff& u \equiv  v ~ (\mod n)\mbox{~and there exist
$w$, $c$ in $\ZZ$ such that $(w,n)=1$}\\
                        &~&   \mbox{and such that $ wa \equiv (v-1)c + b ~(\mod n)  $}.
                        \end{eqnarray*}
\end{fact}

Our objective is to determine a canonical representation of each
conjugacy class of the group $\BB(\ZZ_n)$ in order to simplify the
task of describing rings of invariants of the polynomial ring
$\ZZ_n[x]$. We  use the following simple theorem  that  is an
extremely special case of Dirichlet's theorem on the infinitude of
primes in arithmetic progressions; see \cite[pages
105--122]{Chandra}.

\bigskip
\begin{theorem} \label{Dirichlet}
If $a,b,n$ are positive integers such that $(a,b) = 1$, then the sequence
$$a + k b : k = 0, 1, 2 \cdots$$
contains an element that is relatively prime with  $n$.
\end{theorem}

\noindent
{\it Proof.}
Let  $r$ be the product of all prime factors of $n$ that do not divide
$b$. Then $(r,b) = 1$. Let $b'$ be an inverse of $b$ mod $r$ and let
$k$ be a non-negative integer such that $k \equiv (1-a)b' ~(\mod ~r)$.
Then $a + kb \equiv 1 ~(\mod r)$. Therefore $(a+kb, r) = 1$. By the
definition of $r$, we conclude that $(a+kb, n) = 1$, as desired.
\hfill $\Box$

\begin{theorem} \label{T30:1}
Let $\alpha, \beta \in \BB(\ZZ_n)$ be given by
\begin{eqnarray*}
\alpha (x) =  ux+a ~,~~~~~\beta (x) =  ux+b
\end{eqnarray*}
where $(u,n) = 1$. If $(a,n) = (b,n)$, then $\alpha$ and $\beta$ are
equivalent. In particular, every element in $\BB(\ZZ_n)$ is
equivalent to one of the form $x \mapsto ux+d$ where $n$ is
divisible by $d$.
\end{theorem}

\noindent {\it Proof.} Let $\sigma \in \BB(\ZZ_n)$ be defined by
$\sigma (x) = ux + d$, where $d =(a,n)$. It suffices to show that
$\alpha$ is equivalent to $\sigma$. Let $a_1 = a/d$ and  $n_1 = n/d$
Since $(a_1,n_1)=1$,  there exists $t$ such that $(a_1+tn_1,d) = 1$.
Also, $(a_1+tn_1,n_1) = (a_1,n_1)=1$. Therefore $(a_1+tn_1,n) = 1$
and hence $v := a_1 + t n_1$ is a unit mod $n$. Also, $vd = a + tn
\equiv a  ~(\mod~n)$. Therefore $\sigma(vx) = v(ux+d) =u(vx) + vd =
u(vx) + a$, and $\sigma \cong \alpha$, as desired. \hfill $\Box$

\bigskip
\begin{theorem} \label{T30:4}
Let $\alpha, \beta \in \BB(\ZZ_n)$ be given by
\begin{eqnarray*}
\alpha (x) =  ux+a ~,~~~~~~\beta (x) =  ux+b
\end{eqnarray*}
where $(u,n) = 1$ and where $n$ is divisible by both $a$ and $b$.
Then $\alpha$ and $\beta$ are equivalent if and only if
$(u-1,a)=(u-1,b)$.
\end{theorem}

\noindent {\it Proof.} If $\alpha$ and $\beta$ are equivalent, then
by Fact \ref{3.0} there exist $w,c$ in $\ZZ$ such that $(w,n) = 1$
and $wa \equiv (u-1)c+b ~(\mod n)$. Since $n$ is divisible by $a$,
it follows that $wa \equiv (u-1)c + b ~(\mod a)$ and $b = ka -
(u-1)c$ for some integer $k$. Thus $b$ is divisible by $(u-1,a)$.
Hence $(u-1,b)$ is divisible by $(u-1,a)$. By symmetry, we conclude
that $(u-1,a) = (u-1,b)$.

Conversely, assume that  $(u-1,a) = (u-1,b) = d$, say. Let $a_1=a/d,
b_1= b/d, r = (u-1)/d$. Then $(a_1,r) = (b_1,r) = 1$. Therefore the
congruence $a_i \xi \equiv b_1 ~(\mod r)$ has a solution $\xi$ that
is necessarily relatively prime with $r$. The sequence $(\xi+kr : k
=0,1,2,\cdots)$ consists of solutions of the given congruence and it
contains infinitely many primes. Therefore, one of these solutions
$w$, say, is a unit mod $n$. Thus there exists $w$ such that $(w,n)
= 1$ and $b_1 = w a_1 + c r$. Multiplying by $d$, we have  $b = w a
+ c (u-1)$. By Fact \ref{3.0},  $\alpha$ and $\beta$ are equivalent.
\hfill $\Box$

\medskip We summarize in Corollary  \ref{T30:5} the conclusions obtained
in Theorems \ref{T30:1} and \ref{T30:4}. In the statement of
Corollary \ref{T30:5} we let
$$
U = \left\{ u \in \{ 1, 2, \cdots, n-1\} ~| ~ (u,n) = 1 \right\}.
$$

\medskip
\begin{corollary} \label{T30:5}
Let $\alpha, \beta \in \BB(\ZZ_n)$ be given by
\begin{eqnarray*}
\alpha (x) =  ux+a ~,~~~~~~\beta (x) =  vx+b
\end{eqnarray*}
where $u$ and $v$ are in $U$.  Then $\alpha$ and $\beta$ are
equivalent if and only if $u=v$ and $(u-1,a,n)=(u-1,b,n)$, where
$(-,-,-)$ is the greatest common divisor of the three numbers.

Consequently, every conjugacy class  in $\BB(\ZZ_n)$ has a unique
representation of the form $x \mapsto ux + a$, where $u \in U$,  and
where  both $u-1$ and $n$ are divisible by  $a$.
\end{corollary}

\medskip
Theorems \ref{T30:6} and \ref{T30:7} yield the explicit formula
given in  Corollary \ref{T30:8} for the number of conjugacy classes
in $\BB (\ZZ_n)$.

\medskip
\begin{theorem} \label{T30:6}
Let $\Psi (n)$ denote the number of conjugacy classes in $\BB
(\ZZ_n)$. Then $\Psi$ is  multiplicative in the sense that $\Psi
(rs) = \Psi (r) \Psi (s)$ for all relatively prime positive integers
$r$ and $s$.
\end{theorem}

\noindent {\it Proof.} If  $r$ and $s$ are  relatively prime, then
the rings $\ZZ_{rs}$ and $\ZZ_{r} \times \ZZ_{s}$ are isomorphic by
the Chinese remainder theorem. By Theorem \ref{products}, the groups
$\BB (\ZZ_{rs})$ and $\BB (\ZZ_{r}) \times \BB (\ZZ_{s})$ are
isomorphic. Denoting the number of conjugacy classes of a group $H$
by $\mu (H)$ and using the fact that $\mu (H \times K) = \mu (H)
~\mu (K)$, we see that
\begin{eqnarray*}
\Psi (rs) &=& \mu (\BB (\ZZ_{rs})) ~=~ \mu (\BB (\ZZ_{r}) \times \BB (\ZZ_{s})) \\
&=&  \mu (\BB (\ZZ_{r}))~ ~ \mu( \BB (\ZZ_{s})) ~=~ \Psi(r)~\Psi(s),
\end{eqnarray*}
as desired. \hfill $\Box$

\medskip
\begin{theorem} \label{T30:7}
Let $\Psi (n)$ denote the number of conjugacy classes in $\BB (\ZZ_n)$,
and let $p$ be a prime. Then
$$\Psi (p^e) = \frac{p^{e-1} - 1}{p -1} + p^e.$$
\end{theorem}

\noindent {\it Proof.} According to Corollary  \ref{T30:5}, $\Psi
(p^e)$ is the number of ordered pairs $(u,a)$, where
$$1 \le u < p^e~,~~(u,p^e) = 1~,~~a | p^e~,~~a | (u-1).$$
Let $S$ be the set of pairs that satisfy these conditions, and
let $S_k$, $0 \le k \le p^e$ be those pairs $(u,a)$ in $S$ for
which $a = p^k$. Then $S = \cup_{k=0}^e S_k$.  Also, it is clear that
if $k \ge 1$, then
\begin{eqnarray*}
(u,a) \in S_k &\ifff&  a = p^k  \mbox{~and~} u = 1 + r p^k \mbox{~where~} r = 0,1,\cdots,p^{e-k}-1.
\end{eqnarray*}
Thus card ($S_k$) = $p^{e-k}$ if $k \ge 1$. Also
\begin{eqnarray*}
(u,a) \in S_0 &\ifff&  a = 1  \mbox{~and~} u \mbox{~is a unit mod $p^e$ in $\{1,2,\cdots,p^e\}$}.
\end{eqnarray*}
Thus card ($S_0$) = $\phi (p^e) = p^{e} - p^{e-1}$.
Therefore
\begin{eqnarray*}
\Psi (p^e) &=& \card(S) ~=~ \card(S_0) + \card(S_1) + \cdots + \card(S_e)\\
&=&     \left(p^{e} - p^{e-1}\right) + p^{e-1} + p^{e-2} + \cdots + 1\\
&=& \frac{p^{e-1} - 1}{p -1} + p^e,
\end{eqnarray*}
as desired. \hfill $\Box$

\begin{corollary} \label{T30:8} Let $n = p_1^{e_1} \cdots p_k^{e_k}$
be the factorization of $n$ as a product of distinct prime powers.
The number of conjugacy classes in $\BB (\ZZ_n)$ is
$$
\Psi (n) = \Psi(p_1^{e_1}) \cdots \Psi(p_k^{e_k}),
$$
where $$\Psi(p_i^{e_i})  = \frac{p_i^{e_i-1} - 1}{p_i -1} +
p_i^{e_i},$$ for each $i$ with $1 \le i \le k$.
\end{corollary}

\begin{example} {\em
The group $\BB(\ZZ_9)$ has order 54 and by Theorem \ref{T30:7},
$\BB(\ZZ_9)$ has $\frac{3-1}{3-1} + 3^2 = 10$ conjugacy classes.
Representatives for these conjugacy classes are
\begin{itemize}
\item $x \mapsto x + 9$, the identity element.
\item $x \mapsto x + 3$, with $x \mapsto x + 6$ as conjugate, so a conjugacy class with 2 elements.
\item $x \mapsto x + 1$, with $x \mapsto x + u$, $u \in
\{2,4,5,7,8\}$,  as conjugates, so a conjugacy class with 6
elements.
\item $x \mapsto 2x + 1$, a conjugacy class with 9 elements.
\item $x \mapsto 4x + 3$,  with $x \mapsto 4x + 6$ and $x \mapsto 4x
+ 9$ as conjugates, so a conjugacy class with 3 elements.
\item $x \mapsto 4x + 1$, a conjugacy class with 6 elements.
\item $x \mapsto 5x + 1$, a conjugacy class with 9 elements.
\item $x \mapsto 7x + 3$, a conjugacy class with 3 elements.
\item $x \mapsto 7x + 1$, a conjugacy class with 6 elements.
\item $x \mapsto 8x + 1$, a conjugacy class with 9 elements.

\end{itemize}

}
\end{example}

\end{document}

\end{document}

\begin{center}
{\large\bf Title}\\
On the number of inequivalent $\ZZ_n$-Automorphisms of $\ZZ_n[x]$
\end{center}
\bigskip
\bigskip

\begin{center}
{\large\bf Authors}\\
\end{center}

$$
\left.
\begin{array}{cc}
\mbox{Jebrel M. Habeb} ~~~~~&~~~~~ \mbox{Mowaffaq Hajja}\\
\mbox{Department of Mathematics}~~~~~&~~~~~\mbox{Department of Mathematics}\\
\mbox{Yarmouk University}~~~~~&~~~~~ \mbox{Yarmouk University}\\
\mbox{Irbid  --  Jordan}~~~~~&~~~~~ \mbox{Irbid  --  Jordan} \\
\mbox{ {\it jhabeb@yu.edu.jo}}~~~~~&~~~~~ \mbox{{\it mhajja@yu.edu.jo}}
\end{array}
\right.
$$
\bigskip
\bigskip
\bigskip

\begin{center}
{\bf \large 2000 Mathematics Subject Classification. 13A50.}
\end{center}

\newpage

\newtheorem{acknowledgement}[theorem]{Acknowledgement}
\newtheorem{algorithm}[theorem]{Algorithm}
\newtheorem{axiom}[theorem]{Axiom}
\newtheorem{case}[theorem]{Case}
\newtheorem{claim}[theorem]{Claim}
\newtheorem{conclusion}[theorem]{Conclusion}
\newtheorem{condition}[theorem]{Condition}
\newtheorem{conjecture}[theorem]{Conjecture}
\newtheorem{corollary}[theorem]{\sc Corollary}
\newtheorem{criterion}[theorem]{Criterion}
\newtheorem{definition}[theorem]{\sc Definition}
\newtheorem{example}[theorem]{\sc Example}
\newtheorem{exercise}[theorem]{Exercise}
\newtheorem{lemma}[theorem]{\sc Lemma}
\newtheorem{notation}[theorem]{Notation}
\newtheorem{problem}[theorem]{Problem}
\newtheorem{proposition}[theorem]{Proposition}

\newtheorem{solution}[theorem]{Solution}
\newtheorem{summary}[theorem]{Summary}

\noindent {\it Proof.}

*********************

\noindent \begin{remark} \label{4.6}  {\em
According to Theorem \ref{4.5}, the fixed subring of $R[x]$ under a
subgroup of $G(R)$ is one of the following rings
    $$\mbox{(i)}~ R[x], ~\mbox{(ii)}~ R[x^2,2x], ~\mbox{(iii)} ~  R[y^2,2y]~
    \mbox{(iv)} ~ R[y+2xg] \mbox{~(where~)} g \in R[y].$$
    The first four are the fixed rings of subgroups of $\BB_x$. In fact, letting
    $$\beta : x \mapsto -x+1~,~~\theta : x \mapsto x+1~,~~\gamma : x \mapsto -x,$$
we easily see that
    \begin{eqnarray*}
    R[x] &=& R[x]^{\{e\}}\\
    R[x^2,2x] &=& R[x]^{\langle \gamma \rangle}\\
    R[y] &=& R[x]^{\langle \beta \rangle}\\
    R[y^2,2y] &=& R[x]^{\langle \theta \rangle}
\end{eqnarray*}